\documentclass[11pt]{amsart}
\usepackage{a4wide}
\usepackage{color}
\usepackage{amsmath}
\usepackage{amsthm}
\usepackage{amsfonts}
\usepackage{amssymb}
\usepackage{xspace,xcolor,enumerate}
\usepackage{graphicx}
\usepackage{fancyhdr, psfrag, setspace, supertabular}
\newenvironment{proofof}[1]{\begin{trivlist} \item {\it Proof
#1:~~}}
  {\qed\end{trivlist}}

\mathchardef\mhyphen="2D
\newcommand\R{{\mathbf{R}}}

\newcommand\primes{{\mathcal{P}}}

\newcommand\F{\mathbf{F}}
\newcommand\Z{{\mathbf{Z}}}
\newcommand\N{{\mathbf{N}}}

\newcommand\ZN{\mathbf{Z}_{N}}

\newcommand\E{\mathbf{E}}

\newcommand\vm{\vec{m}}
\newcommand\vh{\vec{h}}

\newcommand\Qf{\mathcal{Q}}
\newcommand\vdc{\mathrm{vdC}}
\newcommand\lcm{\mathrm{lcm}}

\theoremstyle{plain}
  \newtheorem{theorem}{Theorem}

  \newtheorem{proposition}{Proposition}
  
  \newtheorem{lemma}{Lemma}
  \newtheorem{claim}{Claim}

\theoremstyle{remark}
  \newtheorem{remark}[subsection]{Remark}

\theoremstyle{definition}
  \newtheorem{definition}{Definition}

\title{Polynomial configurations in the primes}

\author{Th\'ai Ho\`ang L\^e}
\address{T. H. L\^e, Department of Mathematics,
The University of Texas at Austin,
1 University Station, C1200,
Austin, TX 78712, USA}
\email{leth@math.utexas.edu}

\author{Julia Wolf}
\address{J. Wolf, \'Ecole Polytechnique,
Centre de Math\'ematiques Laurent Schwartz,
91128 Palaiseau, France}
\email{julia.wolf@cantab.net}

\begin{document}

\begin{abstract}
The Bergelson-Leibman theorem states that if $P_1, \dots, P_k \in \Z[x]$, then any subset of the integers of positive upper density contains a polynomial configuration $x+P_1(m), \dots, x+P_k(m)$, where $x,m \in \Z$. Various generalizations of this theorem are known. Wooley and Ziegler showed that the variable $m$ can in fact be taken to be a prime minus 1, and Tao and Ziegler showed that the Bergelson-Leibman theorem holds for subsets of the primes of positive relative upper density. Here we prove a hybrid of the latter two results, namely that the step $m$ in the Tao-Ziegler theorem can be restricted to the set of primes minus 1.
\end{abstract}

\maketitle

%:sec:intro
\section{Introduction}\label{sec:intro}

Roughly twenty years after the ergodic theoretic proof of Szemer\'edi's theorem on long arithmetic progressions in dense subsets of the integers by Furstenberg \cite{furstenberg}, Bergelson and Leibman \cite{bl} proved the following celebrated polynomial generalization.

%:thm:bl
\begin{theorem}[Bergelson-Leibman]\label{thm:bl}
Let $P_1, \ldots, P_{k}$ be polynomials in $\Z[x]$ such that $P_i(0)=0$ for $i=1, \ldots, k$. Then any subset of the integers of positive relative upper density contains a configuration of the form $a+P_1(d),\ldots, a+P_{k}(d)$, where $a,d$ are integers, $d \neq 0$.
\end{theorem}

More recently, Tao and Ziegler \cite{tao-ziegler} proved Theorem \ref{thm:bl} for dense subsets of the primes, using the general transference strategy of Green and Tao \cite{gt-primes}.

%:thm:tz
\begin{theorem}[Tao-Ziegler]\label{thm:tz}
Let $P_1, \ldots, P_{k}$ be polynomials in $\Z[x]$ such that $P_i(0)=0$ for $i=1, \ldots, k$. Then any subset of the primes of positive relative upper density contains
a configuration of the form $a+P_1(a),\ldots, a+P_{k}(a)$, where $a,d$ are integers, $d \neq 0$.
\end{theorem}
Here for any subset $A$ of the set of primes $\mathcal{P}$, the relative upper density $\overline{d}_{\mathcal{P}}(A)$ of $A$ in $\mathcal{P}$ is defined as
\[ \overline{d}_{\mathcal{P}}(A)=\overline{\lim}_{N \rightarrow \infty} \frac{|A\cap [N]|}{|\mathcal{P} \cap [N]|}\]
In a recent preprint, Wooley and Ziegler \cite{wooley-ziegler} showed that the step $d$ of the polynomial progression in Theorem \ref{thm:bl} can be taken to be a shifted prime.

%:thm:wz
\begin{theorem}[Wooley-Ziegler]\label{thm:wz}
Let $P_1, \ldots, P_{k}$ be polynomials in $\Z[x]$ such that $P_i(0)=0$ for $i=1, \ldots, k$. Then any subset of the integers of positive relative upper density contains a configuration of the form $a+P_1(p-1),\ldots, a+P_{k}(p-1)$, where $a$ is an integer and $p$ is prime. The same is true if we replace $p-1$ with $p+1$.
\end{theorem}

A generalization to polynomials in several variables, with a simpler proof, was subsequently obtained by Frantzikinakis, Host and Kra \cite{fhk}.

Our goal in this paper is to establish the following hybrid of Theorems \ref{thm:tz} and \ref{thm:wz}. 

%:thm:main
\begin{theorem}\label{thm:main}
Let $P_1, \ldots, P_{k}$ be polynomials in $\Z[x]$ such that $P_i(0)=0$ for $i=1, \ldots, k$. Then any subset of the primes of positive relative upper density contains a configuration of the form $a+P_1(p-1),\ldots, a+P_{k}(p-1)$, where $a$ is an integer and $p$ is prime. The same is true if we replace $p-1$ with $p+1$.
\end{theorem}
In other words, we claim that the step $d$ in Theorem \ref{thm:tz} can be restricted to be of the form $p-1$ (or $p+1$). In fact, our proof shows that there are infinitely many such configurations, and we are able to give a lower bound on their number which is of the order of magnitude predicted by the Bateman-Horn conjecture. Previously, the question about the existence of such configurations has also been posed as Conjecture 1.2 in \cite{li-pan}.

Our method is very similar to that employed in \cite{fhk}, in the sense that we compare an average over the integers to an average along the shifted primes using multiple applications of van der Corput's lemma and a PET induction scheme. However, we proceed quantitatively in the spirit of \cite{tao-ziegler}, and rely on a refined analysis of the correlation properties of the pseudorandom measure from that paper.

The rest of the article is structured as follows. In Section \ref{sec:prelim} we set up our notation, and in Section \ref{sec:overview} we reduce Theorem \ref{thm:main} to the more technical Propositions \ref{prop:gvn} and \ref{prop:extra}. We study a simple example of Proposition \ref{prop:gvn} in Section \ref{sec:toy}, and follow it up in Section \ref{sec:polyforms} with a discussion of a modified polynomial forms condition that arises from the example, together with an outline of the proof of Proposition \ref{prop:extra}. The technical details of the proof of Proposition \ref{prop:extra} can be found in an appendix. Finally, the general case of Proposition \ref{prop:gvn} is proved in Section \ref{sec:gencase} using the now standard PET induction scheme. 

\textbf{Acknowledgements.} Work on this project began during the first author's visit to \'Ecole polytechnique, and he would like to thank the Centre de Math\'ematiques Laurent Schwartz for its hospitality. The authors would also like to thank Terence Tao for helpful discussions.
%:sec:prelim
\section{Preliminaries}\label{sec:prelim}
We assume some familiarity with the work of Green and Tao \cite{gt-primes} and Tao and Ziegler \cite{tao-ziegler}, as well as with the definition and basic properties of the Gowers uniformity norms. Here we only briefly remind the reader of the most important definitions, lemmas and parameter settings from those papers. The experienced reader is encouraged to skip this section and consult it later as the need arises.

Landau's $O,o$ and Vinogradov's $\ll, \gg$ notation are given their usual asymptotic meaning. That is, for two quantities $X,Y$, we write $X \ll Y, Y \gg X$, or $X = O(Y)$ if we have a bound $|X| \leq CY$ for some constant $C$. If $C$ depends on other parameters such as $k$, then this dependence is indicated as $X \ll_k Y, Y \gg_k X$, or $X=O_k(Y)$. By $o(1)$ we denote a quantity that goes to 0 as $N \rightarrow \infty$. If this quantity depends on other parameters such as $k$, then this dependence is sometimes indicated as $o_k(1)$.

Throughout the paper, we fix a system of polynomials $P_1,\ldots, P_k$ with integer coefficients, each vanishing at 0.  We can certainly assume that these polynomials are distinct. We also fix a subset $A  \subseteq \mathcal{P}$ satisfying $\overline{d}_{\mathcal{P}}(A)=\delta_0>0$. All implicit constants are allowed to depend on $\delta_0, P_1, \dots, P_k$.

To get around the fact that the primes are not equidistributed with respect to small moduli, we let $w \ll \log \log \log N$ be any sufficiently slowly growing function in $N$, and let $W=\prod_{p < w}p$ be the product of the primes less than $w$, so that $W \ll \log \log N$. Eventually, just as in \cite{gt-primes} and \cite{tao-ziegler}, we will be able to take $w$ be a sufficiently large constant, see the discussion in Section \ref{sec:remarks}. %\Jnote{Changed.}

It follows from the assumption on the density of $A$ that there is an infinite sequence of integers $N'$ going to infinity such that 
\[|A \cap [N']| >\frac{1}{2} \delta_0 \frac{N'}{\log N'}.\]
We set $N=\lfloor N'/2W\rfloor$, and observe that the asymptotic limit as $N\rightarrow \infty$ is equivalent to the asymptotic limit as $N'\rightarrow \infty$. By the pigeonhole principle, we can find $b=b(N) \in [W]$ coprime to $W$ such that 
\[ |\{ x \in [N/2] : Wx+b \in A\}| \gg \frac{W}{\phi(W)} \frac{N}{\log N},\]
where $\phi$ is Euler's totient function.

The expression $\E_{y \in Y}f(y)$ denotes the average of a function $f$ over a finite set $Y$. Borrowing notation from ergodic theory, we also write $\int_{X}$ for $\E_{x \in X}f(x)$ and $Tf(x)=f(x-1)$. For convenience we set $X$ equal to the cyclic group $\ZN$. The fact that the elements $x$ we consider are restricted to lie in the interval $[N/2]$ ensures that there is no problem with wrap-around in $X$.

For a modulus $W$ and a residue $1 \leq b \leq W$ coprime to $W$, let us define 
\begin{equation} \label{eq:Lambda}
\Lambda_{W,b;N}(n)= \left\{ 
                                                                   \begin{array}{ll}
                                                                     \frac{\phi(W)}{W}  \log(Wn+b),& \hbox{if $Wn+b$ is prime and $1 \leq n \leq N$,}\\
                                                                     0 & \hbox{otherwise.}
                                                                   \end{array}
                                                                 \right.
\end{equation}                                                       
For the purposes of this paper, a \textit{measure} is a non-negative function $\nu: X \rightarrow [0, \infty)$ satisfying $\int_X \nu = 1+ o(1)$ and the pointwise bound $\nu = O_{\epsilon}(N^{\epsilon})$ for any $\epsilon>0$. The measures we will be working with are of the form
\begin{equation}\label{eq:nu}
\nu_{W,b}(n)=\frac{\phi(W)}{W} \log R \left( \sum_{m|Wn+b} \mu(m) \chi\left( \frac{\log m}{ \log R} \right)  \right)^2,
\end{equation}
where $\mu$ is the M\"obius function and $\chi$ is an even smooth function supported on $[-1,1]$ satisfying
\[\int_{0}^1 | \chi'(t)|^2 dt =1.\]
In \cite{tao-ziegler}, Tao and Ziegler defined a \emph{pseudorandom measure} to be a measure satisfying two technical conditions known as the \textit{polynomial forms} condition and the \textit{polynomial correlation} condition, and they showed that $\nu_{W,b}$ as defined in (\ref{eq:nu}) satisfies both of these. We refer the reader to the precise definitions of the polynomial forms and the polynomial correlation condition in \cite[Definitions 3.6 and 3.9]{tao-ziegler}. In this paper, we will need a variant of the polynomial forms condition for pairs of pseudorandom measures, which we call the \textit{extra condition}. It will be given in Section \ref{sec:polyforms}, where we also verify that this extra condition is satisfied by a pair $\nu_{W,b_1}, \nu_{W,b_2}$ for potentially distinct $b_1, b_2$.

Let us list the important remaining parameters.
\begin{itemize}
\item Let $d_0=\max_{1 \leq i \leq k} \deg P_i$ denote the maximal degree of the polynomials.
\item Let $M=N^{\eta_0}$ be the ``coarse scale", which serves as a bound for the step of the polynomial progression. We can take $\eta_0$ to be any positive number less than $1/2 d_0$.
\item Let $0<\eta_1 \ll \eta_0$ be a tiny parameter, depending on $P_1, \ldots, P_k$, which controls the degree of pseudorandomness of a measure $\nu$.
\item Let $0<\eta_2 \ll \eta_1 / d_0$, and $R=N^{\eta_2}$ be the sieve level which is used in the construction of $\nu$.
\end{itemize}
We do not explicitly specify the parameters $\eta_1, \eta_2$, but insist that they depend only on the system $P_1, \ldots, P_k$ and are chosen sufficiently small to accommodate all our estimates (notably those arising from the PET induction).
Note that Tao and Ziegler also needed the ``fine scale" $H=N^{\eta_7}$, but we shall not need it here. In this sense our work is much simpler than \cite{tao-ziegler}.

For completeness, we state two basic and well-known inequalities we shall use repeatedly.
%:lem:cs
\begin{lemma}[Cauchy-Schwarz]\label{lem:cs}
Let $A, B$ be sets, let $f,F$ be functions on $A$ and let $g$ be a function on $A\times B$. If $|f|\leq F$ pointwise, then
\[ \left| \E_{a \in A, b \in B}f(a)g(a,b) \right|^2 \leq \E_{a \in A} F(a) \E_{a \in A} F(a) \left|\E_{b \in B}g(a,b) \right|^2.\]
\end{lemma}

%:lem:vdc
\begin{lemma}[van der Corput]\label{lem:vdc}
Let $(x_m)_{m \in \Z}$ be a real-valued sequence satisfying $x_m = 0$ outside the interval $[M]$. Then
 \[\left|\E_{m \in [M]}x_m \right|^2 \ll \E_{|h| < M} \E_{m \in [M]} x_{m}x_{m+h}.\]
\end{lemma}
This lemma follows by simply expanding out the square, and is reminiscent of \cite[Lemma A.1]{tao-ziegler}. Note that, in contrast with \cite[Lemma A.1]{tao-ziegler} where $m$ and $h$ are on different scales, in our situation $h$ and $m$ are on the same scale. This fact is important for us since it will make the Gowers norms appear. 

%:sec:overview
\section{Overview of the proof}\label{sec:overview}

The main result of Tao and Ziegler is the following \cite[Theorem 2.3]{tao-ziegler}. We shall use it as a black box in the sequel, although we will need to delve into the details of the proof in a later part of the argument.

%:thm:upsp
\begin{theorem}[Uniform polynomial Szemer\'edi theorem in the primes]\label{thm:upsp} Let $\nu$ be a pseudorandom measure on $X$. If $f$ is a function on $X$ such that $0\leq f \leq \nu$, $\int_X f \geq \delta$, then 
\[ \E_{m \in [M]} \int_{X} T^{P_1(Wm)/W}f \ldots T^{P_k(Wm)/W}f \geq c(\delta) - o(1)\]
for some constant $c(\delta)>0$ depending on $\delta$. 
\end{theorem}

Just as in \cite{fhk}, we will also need the following deep result from Green and Tao's programme of counting linear patterns in primes (see \cite{gt2, gt3, gtz}, but also \cite[Theorem 2.2]{fhk}).

%:thm:gt
\begin{theorem}[Green-Tao, Green-Tao-Ziegler]\label{thm:gt} For every $d \in \Z^{+}$, we have
\[
\lim_{N \rightarrow \infty} \max_{\substack{1 \leq b < W, \\ (b,W)=1}} \| \Lambda_{W,b;N} - 1_{[N]}\|_{U^d(\Z_{(2d+1)N})}=0.
\]
\end{theorem}

Our main result will be deduced from two statements, the first of which is analogous to \cite[Lemma 3.5]{fhk}. The ``extra condition" mentioned in the hypotheses of Proposition \ref{prop:gvn} below is quite technical, and will be defined in Section \ref{sec:polyforms} (Definition \ref{def:extracon}).

%:prop:gvn
\begin{proposition}\label{prop:gvn} Let $\nu_1,\nu_2$ be a pair of pseudorandom measures on $X$ satisfying the extra condition. Suppose that $f_1, \ldots, f_k$ are functions on $X$ with $|f_i| \leq \nu_1$ for $i=1, \dots, k$, and that $a$ is a weight on $X$ with support in $[M]$ such that $|a|\leq 1+\nu_2$. Then 
\[ \E_{m \in [M]} \int_{X} a(m) T^{P_1(Wm)/W}f_1 \cdots T^{P_k(Wm)/W}f_k = O(\|a\|_{U^{d}(\Z_{(2d+1)M})}) + o(1),\]
where $d$ is an integer depending only on the system of polynomials.
\end{proposition}

We shall also show that a pair of measures satisfying the hypotheses of Proposition \ref{prop:gvn} actually exists.

%:prop:extra
\begin{proposition}\label{prop:extra}
For any $b \neq 0$ coprime to $W$, the pair $\nu_1=\nu_{W,b},\nu_2=\nu_{W,1}$ of pseudorandom measures satisfies the extra condition.
\end{proposition}

\begin{remark}
If we were only interested in configurations inside the full set of primes (rather than subsets of positive relative density), this proposition and the needed extra condition would be slightly simpler to state, and to prove. However, we need to be able to take potentially distinct distinct residue classes for the pseudorandom measures governing $a$ and the $f_i$s since the residue class $b$ mod $W$ on which the set $A$ is dense was chosen by the pigeonhole principle.
\end{remark}

To conclude this section, we show how Theorem \ref{thm:main} follows from Propositions \ref{prop:gvn} and \ref{prop:extra}. 

\begin{proofof}{of Theorem \ref{thm:main} assuming Propositions \ref{prop:gvn} and \ref{prop:extra}} 
Suppose that we are given a subset $A \subseteq \primes$ of relative upper density $\delta_0$. We shall let $f_1= \dots=f_k =f$, where 
\begin{equation}\label{eq:fis}
f(x) = \left\{ 
\begin{array}{ll}
\frac{\phi(W)}{W}  \log R & \hbox{if $R \leq x  \leq N/2$ and $Wx+b\in A$,}\\
0 & \hbox{otherwise.}
\end{array}
\right.
\end{equation}
As remarked in Section \ref{sec:prelim}, by the pigeonhole principle we can choose $b$ such that $\int_{X} f \gg \delta_0$ provided that $N$ is sufficiently large. Set $\nu_1=\nu_{W,b}$ so that $0 \leq f \leq \nu_1$. 
Let 
\[
g(x) = \left\{ 
\begin{array}{ll}
\frac{\phi(W)}{W}  \log (Wx+1) & \hbox{if $R \leq x  \leq M$ and $Wx+1$ is prime,}\\
0 & \hbox{otherwise.}
\end{array}
\right.
\]
(in other words, $g$ is the same function as $\Lambda_{W,1;M}$ except on $[R]$), then there is a constant $\alpha$ such that $0 \leq \alpha g \leq \nu_2=\nu_{W,1}$.

Set $a= \alpha (g - 1_{[M]})$, so that $|a|\leq 1+\nu_2$. Proposition \ref{prop:extra} states that the pair $\nu_1,\nu_2$ satisfies the extra condition.

Applying Proposition \ref{prop:gvn} with these choices yields
\[ \E_{m \in [M]} \int_{X} a(m) T^{P_1(Wm)/W}f \cdots T^{P_k(Wm)/W}f =O( \|g - 1_{[M]} \|_{U^{d}(\Z_{(2d+1)M})}) + o(1).\]
Since $g$ and $\Lambda_{W,1;M}$ differ on a negligible subset of $\Z_{(2d+1)M}$, Theorem \ref{thm:gt} tells us that the right-hand side can be made arbitrarily small if $N$ is sufficiently large.
By Theorem \ref{thm:upsp} we also have
\[\E_{m \in [M]} \int_{X} T^{P_1(Wm)/W}f \cdots T^{P_k(Wm)/W}f \geq c(\delta_0) - o(1).\] 
Since $P_i(Wm)$ is much less than $N/2$ for $m \in [M]$ and the progressions therefore cannot wrap around the group $\Z_N$, we can replace the average over $X$ with the average over [N] and find that
\begin{equation}\label{eq:end}
\E_{m \in [M]} \E_{x \in [N]} g(m) f(x+ P_1(Wm)/W ) \cdots f(x+P_k(Wm)/W) \geq c(\delta_0)-o(1).
\end{equation}
It remains to replace $g$ by a suitable indicator function for the primes congruent to 1 mod $W$. Since we are only looking for a lower bound, this is straightforward to accomplish. Indeed, we see that the left-hand side of (\ref{eq:end}) is bounded above by
\[\frac{1}{MN} \frac{\phi(W)}{W} \log (WM+1) \cdot \left(\frac{\phi(W)}{W} \right)^k(\log R)^k\]
times the number of pairs $(m,x) \in [M] \times [N]$ such that $Wm+1 \in \primes$, $x+P_i(Wm)/W \in [N/2] $ and $Wx+b+P_i(Wm) \in A$ for all $i=1,\dots, k$.
This is equivalent to saying that for sufficiently large $N$, the number of pairs $(p,x) \in [WM+1] \times [N]$ satisfying $p \in \primes, p \equiv 1 (W)$ such that $x+P_i(p-1)/W \in [N/2]$ and $Wx+b+P_i(p-1) \in A$ for all $i=1,\dots, k$ is at least
\[ \left(c(\delta_0)-o(1)\right)\frac{MN}{(\log M) (\log R)^k} \left(\frac{W}{\phi(W)}\right)^{k+1}.\]
But the right-hand side tends to infinity with $N$, concluding the proof of Theorem \ref{thm:main}.
\end{proofof}

It thus suffices to prove Propositions \ref{prop:gvn} and \ref{prop:extra}.

%:sec:toy
\section{A toy example}\label{sec:toy}
In this section we will study the toy example of the configuration $x, x+(p-1)^2$, and use it to motivate the definition of the extra condition in the subsequent section. (Note, however, that the existence of this particular configuration in the primes already follows from the work of Li and Pan \cite{li-pan}. A more general result was recently proved by Rice \cite{rice}.) For simplicity we assume here that $W=1$.

Let $\nu_1,\nu_2$ be a pair of pseudorandom measures. Suppose that $f_0, f_1$ are positive functions satisfying $f_0, f_1\leq \nu_1$ and suppose further that the weight $a$, which is supported on $[M]$, satisfies $|a|\leq \nu_2$. We shall show that, under an additional assumption on $\nu_1,\nu_2$, we can prove the estimate
 \begin{equation}\label{eq:gvn1}
  E=\int_{X} \E_{m \in [M]} a(m) f_0(x)T^{m^2}f_1(x) = O(\|a\|_{U^3(\Z_{7M})})+ o(1).
 \end{equation}
Let us first eliminate $f_0$ from the average $E$. By Cauchy-Schwarz, we have
\begin{equation*}
  E^2 \leq \int_X \nu_1(x) \int_{X} \nu_1(x) \left| \E_{m\in [M]} a(m)T^{m^2}f_1(x) \right|^2.
 \end{equation*}
Recalling that $\int_X \nu_1 = 1+ o(1)$ and using van der Corput, we have
\begin{eqnarray*}
 E^2 &\ll& (1+o(1))\int_{X} \E_{\substack{m\in [M],\\ |h|<M}} a(m)a(m+h) \nu_1 T^{m^2}f_1T^{(m+h)^2}f_1 + o(1)\\
 &=& (1+o(1))\int_{X} \E_{\substack{m\in [M],\\ |h|<M}} a(m)a(m+h) T^{-m^2} \nu_1 f_1 T^{2mh+h^2}f_1 + o(1),
\end{eqnarray*}
where in the second line we shift the variable $x$ by $R_1(m)=-m^2$, thus making the term $f_1$ appear, rather than a shift of $f_1$. Let
\[ E_1=\int_{X} \E_{\substack{m\in [M],\\ |h|<M}} a(m)a(m+h) T^{R_1(m)} \nu_1 f_1 T^{2mh+h^2}f_1.\]
Note that the system of shifts of $f_1$ appearing in $E_1$, namely $0, 2mh+h^2$, is ``simpler'' than the 
system in $E$ in the sense that the polynomials are now linear in $m$.
Next, we want to eliminate the new shift of $f_1$ from $E_1$. Again, by Cauchy-Schwarz, we have
\begin{eqnarray*}
 E_1^2 &\leq& \int_X \nu_1(x) \int_{X} \nu_1(x) \left| \E_{\substack{m\in [M],\\ |h|<M}} a(m)a(m+h) T^{R_1(m)} \nu_1 T^{2mh+h^2}f_1 \right|^2,
\end{eqnarray*}
and by van der Corput this is
\begin{eqnarray*}
&\ll& (1+o(1)) \int_X \nu_1 \E_{\substack{m\in [M],\\ |h|,|k|<M}} a(m)a(m+k)a(m+h)a(m+h+k) \\
& & \qquad \times T^{R_1(m)} \nu_1 T^{R_1(m+k)} \nu_1 T^{2mh+h^2}f_1 T^{2(m+k)h+h^2}f_1 + o(1).
\end{eqnarray*}
Let the last integral be $E_2$. Shifting $x$ by $R_2(m,h)=-(2mh+h^2)$ to make $f_1$ appear, we obtain
\begin{eqnarray*}
E_2 &=& \int_X \E_{\substack{m\in [M],\\ |h|, |k| <M}} a(m)a(m+k)a(m+h)a(m+h+k) \\
& & \qquad\times T^{R_2(m,h)}\nu_1 T^{R_1(m)+R_2(m,h)} \nu_1 T^{R_1(m+k)+R_2(m,h)} \nu_1 f_1 T^{2kh}f_1.
\end{eqnarray*}
Again, the system of shifts of $f_1$ is now simpler than the previous one, in that it does not depend on $m$ at all. 
(Tao and Ziegler deduced from this step their generalized von Neumann inequality, which bounds $E$ in terms of an averaged local Gowers norm of $f_1$.) We repeat the same process one more time to eliminate $f_1$ completely from the average. Indeed, by Cauchy-Schwarz, $E_2^2$ is less than or equal to
\begin{eqnarray*}
& & \Big( \E_{|h|,|k| < M} \int_{X} \nu_1 T^{2kh}\nu_1 \Big) \Big( \E_{|h|,|k| < M} \int_{X} \nu_1 T^{2kh}\nu_1 \Big| \E_{m \in [M]} a(m)a(m+k)a(m+h)a(m+h+k) \\
&&\qquad \times T^{R_2(m,h)}\nu_1 T^{R_1(m)+R_2(m,h)} \nu_1 T^{R_1(m+k)+R_2(m,h)} \nu_1 \Big|^2 \Big) .
\end{eqnarray*}
Since $\nu_1$ satisfies the polynomial forms condition \cite[Definition 3.6]{tao-ziegler}, the first factor is $1+o(1)$. By van der Corput, the second factor is at most
\begin{eqnarray*}
& & \E_{\substack{m\in [M],\\ |h|,|k|,|l|<M}} \prod_{\omega \in \{0,1\}^3}a(m+ \omega \cdot (l,k,h)) \int_{X} T^{R_2(m,h)}\nu_1 T^{R_1(m)+R_2(m,h)} \nu_1 T^{R_1(m+k)+R_2(m,h)}\nu_1\\
& & \qquad \times T^{R_2(m+l,h)}\nu_1 T^{R_1(m+l)+R_2(m+l,h)} \nu_1 T^{R_1(m+k+l)+R_2(m+l,h)} \nu_1.
\end{eqnarray*}
If it were not for the presence of the integral, then this would be equal to 
\[
\E_{\substack{m\in [M],\\ |h|,|k|,|l|<M}} \prod_{\omega \in \{0,1\}^3}a(m+ \omega \cdot (l,k,h)),
\] 
which would give us the desired estimate, since the latter quantity is bounded above by a constant times $\| a \|_{U^3(\Z_{7M})}^8$.  Indeed, we trivially have
\[
\| a \|_{U^3(\Z_{7M})}^8 \gg \frac{1}{M^4} \sum_{m, h, k, l \in \Z_{7M}} \prod_{\omega \in \{0,1\}^3}a(m+ \omega \cdot (l,k,h)).
\]
Identifying $\Z_{7M}$ with the integers in $(-3M,4M]$, we see that since $a$ is supported on $[M]$, the term $\prod_{\omega \in \{0,1\}^3}a(m+ \omega \cdot (l,k,h))$ is non-zero only if $m \in [M]$ and $|h|, |k|,|l|<M$. For these $m, h, k, l$, the representative of $m+ \omega \cdot (l,k,h)$ in $(-3M,4M]$ is $m+ \omega \cdot (l,k,h)$ itself, for any $\omega \in \{0,1\}^3$. 
Thus 
\[
\| a \|_{U^3(\Z_{7M})}^8 \gg \frac{1}{M^4} \sum_{\substack{m\in [M],\\ |h|,|k|,|l|<M}} \prod_{\omega \in \{0,1\}^3}a(m+ \omega \cdot (l,k,h)),
\]
where $m, h, k, l$ are now elements of $\Z$, and the claimed bound follows after renormalization.

To continue, let us write $\vm=(m,h,k,l) \in \Z^4$, the integral as $\int_{X} \prod_{i=1}^6 T^{Q_i(\vm)}\nu_1$, where $Q_i \in \Z[\vm]$ for $i=1, \ldots, 6$,
and the weight in front of the integral as $\prod_{j=1}^8 a(L_j(\vm))$, where for $j=1,\dots, 8$ the $L_j$ are linear forms defining the eight vertices of the parallelepiped in the $U^3$ norm.
Also, let $\Omega_M=\{(m,h,k,l) \in \Z^4: m \in [M], |h|, |k|, |l|<M \}$.
We want to show that
\begin{equation} 
F  =  \E_{\vm \in \Omega_M} \prod_{j=1}^{8}a(L_j(\vm)) \times \left( \int_{X} \prod_{i=1}^6 T^{Q_i(\vm)}\nu_1(x)-1\right) = o(1).
\end{equation}
Recalling that $|a| \leq \nu_2$, by Cauchy-Schwarz we have that
\begin{eqnarray*}
|F|^2  &\ll & \|\nu_2\|_{U^{3}(\Z_{7M})}^8 \E_{\vm \in \Omega_M} \prod_{j=1}^{8} \nu_2(L_j(\vm)) \times \left( \int_{X} \prod_{i=1}^6 T^{Q_i(\vm)}\nu_1(x)-1\right)^2.
\end{eqnarray*}
Since $\nu_2$ is pseudorandom, we have $\|\nu_2\|_{U^{3}(\Z_{7M})}=1+o(1)$. (Note that a priori 
$\nu_2$ is a pseudorandom measure with respect to $N$, but by choosing $R$, the sieve level in the definition of $\nu_2$, sufficiently small, we can ensure that $\nu_2$ is also pseudorandom with respect to $M$.)

By squaring out  $\left( \int_{X} \prod_{i=1}^6 T^{Q_i(\vm)}\nu_1(x)-1\right)^2$, we see that it suffices to show that
\begin{equation} \label{eq:last}
\E_{\vm \in \Omega_M} \prod_{j=1}^{8} \nu_2(L_j(\vm)) \times \left( \int_{X} \prod_{i=1}^6 T^{Q_i(\vm)}\nu_1(x)\right)^k=1+o(1)
\end{equation}
for $k=0,1,2$. But it is precisely expressions of this type that will be governed by our new ``extra condition", which we shall formally introduce in the next section.

\begin{remark} It is well known that if $\nu_2$ is a pseudorandom measure, then so is $(\nu_2+1)/2$. It will be easy to see that if the pair $\nu_1, \nu_2$ satisfies the extra condition, then so does the pair $\nu_1,(\nu_2+1)/2$. The above result therefore also applies to the case where $|a| \leq \nu_2 +1$, which is what we need in the proof of Theorem \ref{thm:main}.
\end{remark}

%:sec:polyforms
\section{A discussion of the polynomial forms condition}\label{sec:polyforms}
Let us recall Tao and Ziegler's definition of the polynomial forms condition \cite[Definition 3.6]{tao-ziegler}, of which the extra condition will be a variant.

%:def:polyforms
\begin{definition}[Polynomial forms condition]\label{def:polyforms}
A measure $\nu: X \rightarrow \R^{+}$ is said to satisfy the \emph{polynomial forms condition} if for any family of polynomials $Q_j \in \Z[m_1, \ldots, m_D]$, $j\in J$, satisfying
\begin{itemize}
 \item the difference $Q_i - Q_{j}$ is not constant for $i \neq j$;
 \item the number of polynomials $|J|$ and the number of variables $l$ are bounded by $1/\eta_1$;
 \item the total degree of each $Q_j$ for $j \in J$ is at most $d_0$, and all coefficients are at most $CW^{d_0}$, where $C$ is a constant (depending on $P_1, \ldots, P_k$);
\end{itemize}
we have
\begin{equation} \label{eq:def1}
\E_{\vh \in \Omega \cap \Z^D} \int_X \prod_{j \in J} T^{Q_j(\vh)}\nu(x) = 1+ o_{\epsilon}(1)
\end{equation}
for any convex body $\Omega \subset \R^D$ of inradius at least $N^{\epsilon}$ and contained in the ball $B(0, M^2)$.
\end{definition}

Inequality (\ref{eq:last}) does not exactly follow from Definition \ref{def:polyforms}, but we still can deduce it using Tao and Ziegler's machinery. 
We make the following general definition.
%:def:extracon
\begin{definition}[Extra condition]\label{def:extracon}
A pair of measures $\nu_1,\nu_2: X \rightarrow \R^{+}$ is said to satisfy the \textit{extra condition} if for any family of polynomials $Q_j \in \Z[m_1, \ldots, m_D]$, $j \in J_1$, and any family of linear forms $L_j: \Z^{D} \rightarrow \Z$, $ j \in J_2$, satisfying
\begin{itemize}
 \item the difference of polynomials $Q_i - Q_{j}$ is not constant for $i \neq j$;
 \item the number of polynomials $|J_1|$ and the number of variables $D$ are bounded by $1/\eta_1$;
 \item the total degree of each $Q_j$ for $j\in J_1$ is at most $d_0$, and all coefficients are at most $CW^{d_0}$, where $C$ is a constant (depending on $P_1, \ldots, P_k$);
 \item the linear forms $L_j$, $j\in J_2$ are pairwise linearly independent;
 \item the number of linear forms $|J_2|$ is bounded by $1/\eta_1$;
 \item the coefficients of each $L_j$ are 0 or 1;
\end{itemize}
we have
 \begin{equation} \label{eq:def2}
  \E_{\vm \in \Omega_{M,D}} \prod_{j \in J_2} \nu_2(L_j(\vm)) \left( \int_{X} \prod_{j \in J_1} \nu_1(x + Q_j(\vm)) \right)^k = 1+o(1)
 \end{equation}
for $k=0,1,2$, where $\Omega_{M,D} = \{(m, h_1, \ldots, h_{D-1}) \in \Z^D: m \in [M], |h_i|<M \textup{ for any } i=1 , \ldots, D-1 ) \}$.
\end{definition}
\begin{remark}
 The extra condition is tailor-made to suit our needs. One could merge it with the polynomial forms condition to obtain a more general statement, but this appears unnecessary.
\end{remark}
Let us now show that for any $b \neq 0$ coprime to $W$, the pair $\nu_1=\nu_{W,b}$, $\nu_2=\nu_{W,1}$ given by (\ref{eq:nu}) satisfies the extra condition. That is, we shall turn to proving Proposition \ref{prop:extra}. To begin with, let us recall some more definitions from \cite{tao-ziegler}.
%:def:gbt
\begin{definition}[Good, bad and terrible primes]\label{def:gbt}
Let $P_j \in \Z[x_1, \ldots, x_D]$, $j \in J$, be a family of polynomials. We say a prime $p$ is \emph{good} with respect to the family $P_j, j \in J$, if
\begin{itemize}
\item the polynomials $P_j \pmod p$, $j \in J$, (considered as elements of $\F_p[x_1, \ldots, x_D]$) are pairwise coprime;
\item for each $j \in J$, there is a variable $x_i$ such that $P_j$ can be expressed as $P_j=P_{j,1}x_i+ P_{j,0}$ where $P_{j,1}, P_{j,0} \in \F_p[x_1, \ldots, x_{i-1},x_{i+1}, \ldots, x_D]$ are such that $P_{j,1}$ 
is non-zero and coprime to $P_{j,0}$.
\end{itemize}
We say $p$ is \emph{bad} if it is not good. We say $p$ is \emph{terrible} if at least one of the $P_j$ vanishes identically mod $p$.
\end{definition}
We shall need the following slight variant of the basic correlation estimate \cite[Proposition 10.1]{tao-ziegler}, in which we now have a pair of pseudorandom measures. 

%:prop:corr
\begin{proposition}[Correlation estimate]\label{prop:corr} 
Write $\nu_1=\nu_{W,b_1}$ and $\nu_2=\nu_{W,b_2}$ with $b_1, b_2 \neq 0$ and coprime to $W$. Let $J_1,J_2 \subset \N$ be two disjoint indexing sets, and let $J=J_1 \cup J_2$. For $j \in J$, let $P_j \in \Z[x_1, \ldots, x_D]$ have degree at most $d$. Let $\Omega$ be a convex body in $\R^D$ of inradius at least $R^{4|J|+1}$. 
Let $\mathcal{P}_b$  be the set of primes $w \leq p \leq R^{\log R}$ which are bad with respect to $((WP_j+b_1)_{j \in J_1},(WP_j+b_2)_{j \in J_2})$, and suppose that there are no terrible primes in the same range. Then
\begin{equation}
\E_{x \in \Omega \cap \Z^{D}} \prod_{j \in J_1} \nu_1(P_j(x))\prod_{j \in J_2} \nu_2(P_j(x)) = 1 + o_{D,J,d}(1) + O_{D,J,d} \left( \textup{Exp} \left( O_{D,J,d} \left( \sum_{p \in \mathcal{P}_b} \frac{1}{p} \right) \right) \right).
\end{equation}
\end{proposition}
Here we have written $\textup{Exp}(x)=\max(e^x-1,0)$, so that $\textup{Exp}(x)\ll x$ when $0 \leq x \ll 1$.

The proof of \cite[Proposition 10.1]{tao-ziegler} generalizes readily to yield Proposition \ref{prop:corr}. However, since it is relatively complex and buried in various appendices of a long paper, we give the details for the convenience of the reader in the appendix to this paper.

To conclude, let us deduce the extra condition from the above correlation estimate.
%:proofof:prop:extra
\begin{proofof}{of Proposition \ref{prop:extra} assuming Proposition \ref{prop:corr}}
Recall that we want to show that
\begin{equation} \label{eq:verify}
\E_{\vm \in \Omega_{M,D}} \prod_{j\in J_2} \nu_2(L_j(\vm)) \left( \int_{X} \prod_{j \in J_1} \nu_1(x + Q_j(\vm)) \right)^k = 1+o(1)
\end{equation}
for $k=0,1,2$, where $\nu_1,\nu_2$ are defined as in Proposition \ref{prop:corr}. If $k=0$, then the integral disappears, and we are left to show that
\[\E_{\vm \in \Omega_{M,D}} \prod_{j \in J_2} \nu_2(L_2(\vm)) = 1+o(1).\]
If we choose $R$ sufficiently small in terms of $M$, then $\nu_2$ is pseudorandom with respect to $M$, and (\ref{eq:verify}) simply follows from Green and Tao's linear forms condition in \cite[Definition 3.1]{gt-primes} in this case. 

Let us now discuss the case $k=2$ (the case $k=1$ is even simpler). First we make a reduction, replacing the integral on $X$ with the average $\E_{x \in [N]}$, thus regarding $\nu_1$ as a function on $\Z$ rather than $X$. 
The values of  $\nu_1(x + Q_j(\vm))$ may be different when $\nu_1$ is regarded as a function on $\Z$ because of the wrap-around effect, but they must agree whenever
$1 \leq x \leq N- O((WM)^{d_0})$.
Recall that we also have the bound $\nu_1 \ll_{\epsilon} N^{\epsilon}$ for any $\epsilon>0$. Thus 
\begin{equation} \label{eq:d}
 \int_{X} \prod_{j \in J_1} \nu_1(x + Q_j(\vm)) -  \E_{x \in [N]} \prod_{j \in J_1} \nu_1(x + Q_j(\vm)) \ll_{\epsilon} (WM)^{d_0} N^{\epsilon-1} ,
\end{equation}
for any $\epsilon>0$, and (\ref{eq:verify}) follows if we can show that
\begin{equation} \label{eq:verify2}
\E_{\vm \in \Omega_{M,D}} \prod_{j \in J_2} \nu_2(L_j(\vm)) \left( \E_{x \in [N]} \prod_{j \in J_1} \nu_1(x + Q_j(\vm)) \right)^2 = 1+o(1).
\end{equation}
Expanding out (\ref{eq:verify2}), we see that it is equivalent to
\begin{equation} \label{expand}
 \E_{\vm \in \Omega_{M,D}, x,x' \in [N]} \prod_{j \in J_2} \nu_2(L_j(\vm)) \prod_{j \in J_1}\nu_1(x + Q_j(\vm)) \nu_1(x' + Q_j(\vm))=1+o(1).
\end{equation}
Now this expression falls within the scope of Proposition \ref{prop:corr}: the polynomials in question are $L_i(\vm), x+Q_j(\vm), x'+Q_j(\vm)$, in variables $\vm, x, x'$. For this system, there is no terrible prime greater than $w$, and the only bad primes greater than $w$ are those dividing $Q_{j}-Q_{j'}$ for some $j \neq j'$.
Therefore, the left hand side of (\ref{expand}) equals $1+o(1)+ O \left( \textup{Exp} \left( O \left( \sum_{p \in \mathcal{P}_b} p^{-1} \right) \right) \right)$, where $\mathcal{P}_b$ denotes the set of primes dividing $Q_{j}-Q_{j'}$ for some $j \neq j'$. 

But just as in the proof of \cite[Corollary 11.2]{tao-ziegler}, if a prime $p$ divides $Q_{j}-Q_{j'}$ for some $j \neq j'$, then $p$ must divide a non-zero difference of the coefficients of the $Q_j$s (recall that $Q_{j}-Q_{j'}$ are not constant for any $j \neq j$). These coefficients are bounded by $O(W^{d_0})$, so that the total product of such $p$ is at most $O(W^{O(1)})$. As a result, the number of $p \in \mathcal{P}_b$ (which are greater than $w$) is at most $\log (O(W^{O(1)}))/\log w = o(\log W)$. But then $\sum_{p \in \mathcal{P}_b} p^{-1} < \sum_{p \in \mathcal{P}_b} w^{-1} =o(1)$, since $\log W \ll w$. 
\end{proofof}

\begin{remark}
The extra condition is in a sense simpler than Tao and Ziegler's full polynomial forms condition, in that all the variables are at scale $M$, whereas the polynomial forms condition (\cite[Theorem 11.1]{tao-ziegler}) makes a statement about convex bodies of inradius as small as $N^{\epsilon}$. This explains why Tao and Ziegler had to do some extra work to prove the polynomial forms condition from the correlation estimate, while for us it is almost immediate.
\end{remark}

%:sec:gencase
\section{The general case}\label{sec:gencase}

As is to be expected, we proceed by PET induction to prove Proposition \ref{prop:gvn} in the general case. We follow the notation in \cite{fhk} with the simplification that the dimension is equal to 1.

Given a family of polynomials $\Qf=(q_1, \dots, q_k)$ in a variable $n$ (and possibly in other variables), the maximum of the degrees of the $q_i$ (with respect to $n$) is called the \emph{degree of the family} $\Qf$. We work with families of polynomials whose degree is smaller than or equal to a fixed number $s$.

For $i=1, \dots, k$, define $\Qf'$ to be the possibly empty set
\[ \Qf'= \{ q_i \in \Qf : q_{i} \textrm{ is constant in $n$ } \}.\] 

Two polynomials are said to be \emph{equivalent} if they have the same degree and the same leading coefficient (in $n$). For $j=1, \dots, s$, let $w_j$ denote the number of distinct non-equivalent classes of polynomials of degree $j$ in $\Qf \setminus \Qf'$. Finally, define the \emph{type} of the family $\Qf$ to be the vector $(w_1, \dots, w_s)$. A family is said to be of type zero if all the $w_j$ are zero, in which case all the polynomials are constant (in $n$). The set of types can be ordered lexicographically, by stipulating that $w=(w_1, \dots, w_s)<w'=(w_1', \dots, w_s')$ if there exists $d$ such that $w_d<w_d'$ and $w_{j}=w_{j}'$ for all $j>d$.

It follows that any decreasing sequence of types of families of polynomials is eventually stationary, and thus any inductive process that reduces the type must eventually stop.

Given a family $\Qf=(q_1, \dots, q_k)$, $q \in \Z[t]$ and $h \in \N$, define following \cite{cfh} the \emph{van der Corput operation} $(q,h) \mhyphen \vdc(\Qf)$ by setting
\[(q,h) \mhyphen \vdc(\Qf)=(S_h\Qf -q,\Qf-q),\]
where $S_hq(n)=q(n+h)$, $S_h\Qf=(S_h q_1, \dots, S_h q_k)$ and $\Qf-q= (q_1-q, \dots, q_k-q)$. 

The crucial observation is the following \cite{bl,fhk}.

%:lem:pet
\begin{lemma}\label{lem:pet}
Let $\Qf$ be a family of polynomials of non-zero type. Then there exists $q \in \Qf \setminus \Qf'$ such that for all $h \in \N$, the family 
$(q,h) \mhyphen \vdc(\Qf \setminus \Qf')$ has strictly smaller type than $\Qf$.
\end{lemma}

For example, in the toy example in Section \ref{sec:toy}, we passed from type $(1,0,1)$ to $(1,1,0)$ to $(2,0,0)$.

In the van der Corput operation, we focus on a single variable $n$. However, we also need to keep in mind that the polynomials in $\Qf$ and $(q,h) \mhyphen \vdc(\Qf)$ are \textit{multivariate}, with a new variable being introduced at each step of the van der Corput operation. This is important when verifying the hypotheses of the (polynomial forms or extra) condition as we apply them to the averages arising throughout.

%:proofof:prop:gvn
\begin{proofof}{of Proposition \ref{prop:gvn}}
Let us recall that we have functions $f_i, i=1,\dots, k$ satisfying $|f_i| \leq \nu_1$, and a weight $a$ supported on $[M]$ satisfying $|a| \leq \nu_2$. We start with the average
\[ E=\E_{m \in [M]} \int_X a(m) \prod_{q_i \in \Qf} T^{q_i} f_i,\]
where the family $\Qf$ initially consists of the polynomials $q_i(m)=P_i(Wm)/W$, $j=1,\dots,k$. Define $\Qf' \subseteq \Qf$ to be the subset of polynomials which are constant in the variable $m$. (By hypothesis on the $P_i$, this actually means that the $q \in \Qf'$ are identically 0, but later on, with additional variables, this is not necessarily the case.)

We first apply the Cauchy-Schwarz inequality to obtain
\begin{equation*}
|E|^2\leq\left( \int_X \prod_{q \in \Qf'} T^q \nu_1\right)\left( \int_X \prod_{q \in \Qf'} T^q \nu_1 |\E_{m \in [M]} a(m) \prod_{q_i \in \Qf\setminus\Qf'} T^{q_i} f_i|^2 \right).
\end{equation*}
By the properties of $\nu_1$, the first integral is $1+o(1)$ (if $\Qf' = \emptyset$, we interpret the empty product as equal to 1), and we can bound the second factor by van der Corput by a constant times
\begin{equation*}
\int_X \prod_{q \in \Qf'} T^q \nu_1 \E_{\substack{m \in [M], \\ |h_1|<M}} a(m)a(m+h_1) \prod_{q_i \in \Qf\setminus\Qf'} T^{q_i(m)} f_i T^{q_i(m+h_1)} f_i +o(1)
\end{equation*}
Now by Lemma \ref{lem:pet}, there exists $q_1 \in \Qf\setminus \Qf'$ such that for all $h_1$, the family $\Qf_1:=(q^{(1)},h_1) \mhyphen \vdc(\Qf\setminus \Qf')$ has strictly smaller type than $\Qf$. We also write $\Qf_1^\dag:=\Qf'-q^{(1)}$ for the recently deceased nodes. With this notation, shifting by $q^{(1)}$ gives, up to an error term, the expression
\begin{equation*}
E_1=\int_X \E_{\substack{m \in [M], \\ |h_1|<M}} a(m)a(m+h_1) \prod_{q \in \Qf_1^\dag} T^q \nu_1 \prod_{q_i \in \Qf_1} T^{q_i}f_{j_i},
\end{equation*}
where the $f_{j_i}$ belong to the set $\{f_1, \dots, f_k\}$. There is no need to keep track of them, so we shall simply write $f$ with no subscript in the sequel.

Write $\Qf_1'$ for those polynomials in $\Qf_1$ of degree 0 in $m$. By Cauchy-Schwarz,
\[ |E_1|^2 \leq \left(\int_X \prod_{q \in \Qf_1'} T^q \nu_1 \right)\left(\int_X \prod_{q \in \Qf_1'} T^q \nu_1 \left|\E_{m,h_1 \in [M]} a(m)a(m+h_1) \prod_{q \in \Qf_1^\dag} T^q \nu_1 \prod_{q \in \Qf_1\setminus \Qf_1'} T^{q}f\right|^2 \right).\]
By the polynomial forms condition, the first factor is $1+o(1)$, and the second can be bounded above by van der Corput as a constant times
\begin{align*}
\int_X \E_{\substack{m \in [M], \\ |h_1|, |h_2|<M}}&a(m)a(m+h_1)a(m+h_2)a(m+h_1+h_2)\\
&\qquad \prod_{q \in \Qf_1^\dag} T^{q(m)} \nu_1 T^{q(m+h_2)}\nu_1 \prod_{q \in \Qf_1\setminus \Qf_1'} T^{q(m)}f T^{q(m+h_2)}f.
\end{align*}
Of course the dependence of the polynomials on $h_1$ is suppressed here. By Lemma \ref{lem:pet}, there is $q^{(2)} \in \Qf_1 \setminus \Qf_1'$ such that the family $\Qf_2 := (q^{(2)},h_2)\mhyphen\vdc(\Qf_1 \setminus \Qf_1')$ has strictly smaller type. Define also $\Qf_2^\dag:=(\Qf_1'-q^{(2)}) \cup (\Qf_1^\dag-q^{(2)})$, leading after rearranging to
\[ E_2= \int_X\E_{\substack{m \in [M], \\ |h_1|, |h_2| <M}}a(m)a(m+h_1)a(m+h_2)a(m+h_1+h_2) \prod_{q \in \Qf_2^\dag} T^q \nu_1 \prod_{q \in \Qf_2} T^{q} f.\]
Continuing in this vein, we set at step $s$ $\Qf_{s+1} := (q^{(s+1)},h_{s+1}) \mhyphen \vdc(\Qf_s \setminus \Qf_s')$ and $\Qf_{s+1}^\dag:=(\Qf_s'-q^{(s+1)}) \cup (\Qf_s^\dag-q^{(s+1)})$, all the while reducing the type. We also set $\Qf_0=\Qf, \Qf^\dag_0=\emptyset$, so that the above recursive definition is valid for all $s \geq 0$. By Lemma \ref{lem:pet} we reach a point, at step $t$ say, where $\Qf_t=\Qf_t'$. In other words, the system is of zero type and all active polynomials are of degree 0 in $m$. 

In order to be able to use the polynomial forms condition at every step, we need to ensure that no two polynomials in $\Qf_s'$ differ by constants. We will in fact prove a slightly stronger statement, which we shall need later.
%:claim:1
\begin{claim} \label{claim:1} For any $0 \leq s \leq t$, no two polynomials in $\Qf_s  \cup \Qf^\dag_{s}$ differ by constants.
\end{claim}
\begin{proofof}{of Claim \ref{claim:1}}
It is easy to see that all our polynomials are 0 when all the variables are 0. Therefore, it suffices to show that all the polynomials in $\Qf_s  \cup \Qf^\dag_{s}$ are distinct as multivariate polynomials. We prove this by induction on $s$. When $s=0$, this follows from our assumption that the polynomials $P_i$ are distinct. Suppose we know already that all the polynomials in $\Qf_s  \cup \Qf^\dag_{s}$ are distinct. By definition, the family $\Qf_{s+1}  \cup \Qf^\dag_{s+1}$ is obtained by subtracting $q^{(s+1)}$ from all the polynomials in the family $S_{h_{s+1}} (\Qf_s \setminus \Qf'_s) \cup (\Qf_s \setminus \Qf'_s) \cup \Qf'_s \cup \Qf_s^{\dag}$. Thus it suffices to show that all polynomials in the latter family are distinct. Recall that the polynomials in $\Qf_s \setminus \Qf'_s$ are distinct and non-constant in $m$. Thus the polynomials in $S_{h_{s+1}} (\Qf_s \setminus \Qf'_s)$ are distinct from each other and from the rest, since they have a new variable, namely $h_{s+1}
 $, and are non-constant in this variable. The remaining polynomials from $(\Qf_s \setminus \Qf'_s) \cup \Qf'_s \cup \Qf_s^{\dag}$ are distinct by induction hypothesis since $(\Qf_s \setminus \Qf'_s) \cup \Qf'_s=\Qf_s$. 
\end{proofof} 

Claim \ref{claim:1} shows that our calculations so far have been valid. Returning to step $t$, we have
\[ E_t=\E_{\substack{m \in [M],\\ |h_1|, \dots, |h_t| < M}} \prod_{ \omega \in \{0,1\}^t} a(m+\omega\cdot (h_1, \dots, h_t))\prod_{q \in \Qf_t^\dag} T^q \nu_1 \prod_{q \in \Qf_t'} T^{q} f,\]
and we apply Cauchy-Schwarz one more time to get
\begin{align*}
|E_t|^2 \leq  & \left(\int_X \E_{|h_1|, \dots, |h_t| < M} \prod_{q \in \Qf_t'} T^{q}\nu_1 \right) \\
& \times \left(\int_X \E_{|h_1|, \dots, |h_t| < M} \prod_{q \in \Qf_t'} T^{q}\nu_1 \left| \E_{m \in [M]}\prod_{ \omega \in \{0,1\}^t} a(m+\omega\cdot (h_1, \dots, h_t))\prod_{q \in \Qf_t^\dag} T^q \nu_1 \right|^2 \right).
\end{align*}
A final van der Corput gives
\[ \E_{\substack{m \in [M],\\ |h_1|, \dots, |h_{t+1}| < M}} \prod_{\omega \in \{0,1\}^{t+1}} a(m+\omega\cdot (h_1, \dots, h_{t+1})) \int_X \prod_{q \in \Qf_t'} T^{q(m)}\nu_1 \prod_{q \in \Qf_t^\dag} T^{q(m)} \nu_1 T^{q(m+h_{t+1})}\nu_1,\]
and we write $F$ for the difference between this expression and $\E_{\substack{m \in [M],\\ |h_1|, \dots, |h_t| < M}} \prod_{\omega \in \{0,1\}^{t+1}} a(m+\omega\cdot (h_1, \dots, h_{t+1}))$. As in the example in Section \ref{sec:toy}, by Cauchy-Schwarz we have that $F^2$ is bounded above by $\|\nu_2\|_{U^{t+1}(\Z_{(2t+3)M})}^{2^{t+1}}$ times the average
\begin{align} \label{eq:genlast}
\E_{m, h_1, \dots, h_{t+1}\in [M]} & \prod_{\omega \in \{0,1\}^{t+1}} \nu_2(m+\omega\cdot (h_1, \dots, h_{t+1}))\\
&\times \left(\int_X \prod_{q \in \Qf_t'} T^{q(m)}\nu_1\prod_{q \in \Qf_t^\dag} T^{q(m)} \nu_1 T^{q(m+h_{t+1})}\nu_1-1 \right)^2.\nonumber
\end{align}
In order to use the extra condition on the last expression, we need to verify that no two polynomials in $\Qf_t^\dag \cup S_{h_{t+1}} \Qf_t^\dag \cup \Qf_t'$ differ by constants. Again, since all these polynomials are 0 when evaluated at 0, it suffices to show that they are distinct. Before seeing this, let us make some observations.
%:claim:2
\begin{claim} \label{claim:2}
For any $0 \leq s \leq t$, for any polynomials $p \in \Qf_s^\dag $ and $q \in \Qf_s \setminus \Qf_s'$, $p-q$ is not constant in $m$.
\end{claim}
\begin{proofof}{of Claim \ref{claim:2}}
For $s=0$ there is nothing to prove. Suppose the claim is true for $s \leq t-1$. Let $p \in \Qf_{s+1}^\dag $ and $q \in \Qf_{s+1} \setminus \Qf_{s+1}'$.
Write $p=u-q^{(s+1)}$ for $u \in \Qf_{s}'\cup \Qf_{s}^\dag$ and $q=v-q^{(s+1)}$ for $v \in (\Qf_{s} \setminus \Qf_{s}') \cup S_{h_{s+1}}(\Qf_{s} \setminus \Qf_{s}') $. It remains to see that $u-v$ is not constant in $m$. If $u \in \Qf_s'$, then $u$ is constant in $m$, but none of the polynomials in $(\Qf_{s} \setminus \Qf_{s}') \cup S_{h_{s+1}}(\Qf_{s} \setminus \Qf_{s}')$ are constant in $m$. Suppose $u \in \Qf_s^\dag$. If $v \in \Qf_s \setminus \Qf_s'$ then $u-v$ is not constant in $m$ by induction hypothesis. If $v \in S_{h_{s+1}}(\Qf_s \setminus \Qf_s')$, we write $v(m)=w(m+h_{s+1})$ for some $w \in \Qf_s \setminus \Qf_s'$. Then $u(m)-w(m+h_{s+1})$
is not constant in $m$, since it is already not constant in $m$ upon setting $h_{s+1}=0$. This proves Claim \ref{claim:2}.
\end{proofof}
%:claim:3
\begin{claim} \label{claim:3}
For any $0 \leq s \leq t$, the polynomials in $\Qf_s^\dag$ are not constant in $m$.
\end{claim}
\begin{proofof}{of Claim \ref{claim:3}}
Indeed, if $s \geq 1$ then by definition we have
\[
 \Qf_s^\dag = (\Qf_{s-1}^\dag-q^{(s)}) \cup (\Qf_{s-1}' - q^{(s)})
\]
for some $q^{(s)} \in \Qf_s\setminus \Qf_s'$. Since the polynomials in $\Qf_{s-1}'$ are constant in $m$, the polynomials in $\Qf_{s-1}' - q^{(s)}$ are not constant in $m$. From Claim \ref{claim:2}, we know that the polynomials in $\Qf_{s-1}^\dag-q^{(s)}$ are not constant in $m$. 
\end{proofof}

From Claim \ref{claim:1} we know that the polynomials in $\Qf_t^\dag \cup \Qf_t'$ are distinct. The polynomials in
$S_{h_{t+1}} \Qf_t^\dag$ have a new variable, namely $h_{t+1}$. From Claim \ref{claim:3} we know that they are not constant in $h_{s+1}$, hence distinct from $\Qf_t^\dag \cup \Qf_t'$.

It is also easy to see that the total degrees of the polynomials appearing in this process are not increased, so they are always at most $d_0$. 
Also, all of their coefficients can be bounded by a constant $C$ times the maximum of the absolute values of the coefficients of the original polynomials (namely $P_i(Wm)/W$), where $C$ depends only on $P_1, \ldots, P_k$. 
It follows that all polynomial expressions in (\ref{eq:genlast}) satisfy the hypotheses of the extra condition. As in Section \ref{sec:toy}, expanding out the expression in (\ref{eq:genlast}) and using the extra condition to see that it equals $o(1)$ concludes the proof of Proposition \ref{prop:gvn}.
\end{proofof}

%:sec:remarks
\section{Concluding remarks}\label{sec:remarks}
Our proof actually gives a lower bound for the number of desired configurations. More precisely, it shows that the number of pairs $(n,p) \in [N]\times [M]$ for which
$n+P_1(p-1), \ldots, n+P_k(p-1)$ and $p$ are all prime is at least $c NM/(\log N)^{k+1}$, as long as $M$ grows like a power of $N$ that is at most $N^{1/2d_0}$. Just as in \cite{tao-ziegler}, this follows from the proof of Theorem \ref{thm:main} since we can choose $w$ to be arbitrarily slowly growing (see also the more detailed discussion at the start of \cite[Section 11]{gt-primes}). This is of the correct order of magnitude if one assumes the Bateman-Horn conjecture. 

One can also see that Proposition \ref{prop:gvn} remains true if we apply it to functions $f_i$ satisfying 
$|f_i| \leq 1$ instead of $|f_i| \leq \nu$. Indeed, under this condition, the proof of Proposition \ref{prop:gvn} is even more straightforward: at each step, in place of the factor $\left( \int_X \prod_{q \in \Qf_i'} T^q \nu_1\right)=1+o(1)$, one simply has a constant $1$. Thus, taking $f_1=\ldots=f_k$ to be the characteristic functions of a set 
$A \subseteq [N]$ of density $\delta$ and using the uniform Bergelson-Leibman theorem \cite[Theorem 3.2]{tao-ziegler}, one actually obtains a slightly different proof of Theorem \ref{thm:wz}. It is, of course, in the same 
spirit as \cite[Theorem 1.2]{fhk}, but has the advantage that it gives a lower bound on the number of configurations. More precisely, we have the following result.

%:prop:qfhk
\begin{proposition}\label{prop:qfhk}
 Let $\kappa>0$. Suppose $N^{\kappa}< M < N^{1/2d_0}$. Then for any $\delta >0$, there is a constant $c(\kappa, \delta)>0$ such that the following holds. 
 Let $A$ be any subset of $[N]$ of density $\delta$. Then $A$ contains at least $c(\kappa, \delta) NM/\log M$ configurations of the form $a+P_1(p-1), \ldots, a+P_k(p-1)$, where $p\leq M$ is a prime.
\end{proposition}

This bound does not follow from \cite{fhk}. On the other hand, the proof of Frantzikinakis, Host and Kra shows that if 
$A \subset \Z$ and $\overline{d}(A)>0$, then 
\[ \overline{d}(A \cap (A-P_1(p-1)) \cdots \cap (A - P_k(p-1))) >0\]
for $p$ in a set of positive relative density in the primes.

%:sec:appendix
\section*{Appendix: The generalized correlation estimate}\label{sec:appendix}

In this appendix we point out the modifications that need to be made to the proof of \cite[Proposition 10.1]{tao-ziegler} to obtain Proposition \ref{prop:corr} (correcting some misprints from \cite{tao-ziegler} in the process).
Recall that we had $\nu_1=\nu_{W,b_1}$ and $\nu_2=\nu_{W,b_2}$ and two disjoint indexing sets $J_1,J_2 \subset \N$, $J=J_1 \cup J_2$. For $j \in J$, we have polynomials $P_j \in \Z[x_1, \ldots, x_D]$ of degree at most $d$. The convex body $\Omega \subset \R^D$ was assumed to have inradius at least $R^{4|J|+1}$. 
We denoted by $\mathcal{P}_b$ the set of primes $w \leq p \leq R^{\log R}$ which are bad with respect to $((WP_j+b_1)_{j \in J_1},(WP_j+b_2)_{j \in J_2})$, and assumed that there are no terrible primes in the same range.

%:proofof:prop:corr
\begin{proofof}{of Proposition \ref{prop:corr}}
We wish to estimate
\[\E_{x \in \Omega \cap \Z^{D}} \prod_{j \in J_1} \nu_1(P_j(x))\prod_{j \in J_2} \nu_2(P_j(x)).\]
Expanding this out in terms of the definitions of $\nu_1, \nu_2$, we find that
\begin{align}
\left(\frac{\phi(W)}{W} \log R\right)^{|J|} \sum_{j\in J} \sum_{m_j, m_j'\geq 1} &\left( \prod_{j \in J} \mu(m_j)\mu(m_j') \chi\left( \frac{\log m_j}{\log R}\right)\chi\left( \frac{\log m_j'}{\log R}\right) \right)\\
& \E_{x \in \Omega \cap \Z^{D}} \prod_{j \in J_1} 1_{\lcm(m_j,m_j') | WP_j(x)+b_1}\prod_{j \in J_2} 1_{\lcm(m_j,m_j') | WP_j(x)+b_2},\nonumber
\end{align}
where $\lcm(a,b)$ denotes the least common multiple of two integers $a$ and $b$. Setting $M=\lcm((m_j)_{j \in J},(m_j')_{j \in J})$, we observe that $M$ can be assumed to be square-free (due to the presence of the M\"obius function) and of size at most $R^{2|J|}$ (due to the restrictions on each $m_j, m_j'$ imposed by the cutoff $\chi$). Each function $x \mapsto 1_{\lcm(m_j,m_j') | WP_j(x)+b_1}$, $x \mapsto 1_{\lcm(m_j,m_j') | WP_j(x)+b_2}$ is periodic with respect to $M \cdot \Z^{D}$, and can therefore be defined on $\Z_M^D$. By \cite[Corollary C.3]{tao-ziegler}, we have
\begin{align*}
 \E_{x \in \Omega \cap \Z^{D}}& \prod_{i=1,2}\prod_{j \in J_i} 1_{\lcm(m_j,m_j') | WP_j(x)+b_i}\\&= \left(1+O\left(R^{-2|J|-1}\right)\right)\E_{y \in \Omega \cap \Z_M^{D}} \prod_{i=1,2}\prod_{j \in J_i} 1_{\lcm(m_j,m_j') | WP_j(y)+b_i},
\end{align*}
where the $O$ error term is easily seen to result in an additive $o(1)$ error, and will therefore be negligible. Setting
\[ \alpha_{(a_j)_{j \in J}}=\E_{y \in \Z_{\lcm((a_j)_{j \in J})}^D} \prod_{i=1,2}\prod_{j \in J_i} 1_{a_j | WP_j(y)+b_i},\]
it therefore suffices to show that
\begin{align}\label{eq:prev}
\left(\frac{\phi(W)}{W} \log R\right)^{|J|} \sum_{j\in J} \sum_{m_j, m_j'\geq 1}&\left( \prod_{j \in J} \mu(m_j)\mu(m_j') \chi\left( \frac{\log m_j}{\log R}\right)\chi\left( \frac{\log m_j'}{\log R}\right) \right)\alpha_{(\lcm(m_j,m_j'))_{j \in J}} \\
&= 1 + o_{D,J,d}(1) + O_{D,J,d} \left( \textup{Exp} \left( O_{D,J,d} \left( \sum_{p \in \mathcal{P}_b} \frac{1}{p} \right) \right) \right).\nonumber
\end{align}
By the Chinese remainder theorem $\alpha$ is multiplicative in the sense that if $\lcm(m_j,m_j')=\prod_p p^{r_j(p)}$, then
\[ \alpha_{(\lcm(m_j,m_j'))_{j \in J}}=\prod_{p} \alpha_{(p^{r_j(p)})_{j\in J}},\]
where the latter is a finite product. Since the $m_j$ are assumed to be squarefree, $r_j(p)$ is either 0 or 1 for each $j$ and each $p$, and we obtain
\[\alpha_{(\lcm(m_j,m_j'))_{j \in J}}=\prod_{p} c_p((WP_j+b_1)_{j \in J_1, r_j(p)=1},(WP_j+b_2)_{j \in J_2, r_j(p)=1}),\]
where the local factor $c_p(P_1, \dots, P_k)$ is defined by
\[ c_p(P_1, \dots, P_k)=\E_{y \in \F_p^D} \prod_{j \in J} 1_{P_j(y) \equiv 0 (p)}.\]
So the left-hand side of (\ref{eq:prev}) becomes
\begin{align}\label{eq:prev2}
\left(\frac{\phi(W)}{W} \log R\right)^{|J|} \sum_{j\in J}& \sum_{m_j, m_j'\geq 1} \left( \prod_{j \in J} \mu(m_j)\mu(m_j') \chi\left( \frac{\log m_j}{\log R}\right)\chi\left( \frac{\log m_j'}{\log R}\right) \right)\\
&\prod_{p\leq R^{\log R}} c_p((WP_j+b_1)_{j \in J_1, r_j(p)=1},(WP_j+b_2)_{j \in J_2, r_j(p)=1})\nonumber,
\end{align}
where we were able to restrict the product to primes less than $R^{\log R}$ because each $m_j$ is bounded by $R$.

We can now replace $\chi$ by terms which are multiplicative in $m_j, m_j'$, using the Fourier expansion
\[ \chi(x)=e^{-x}\int_{-\infty}^{\infty} \phi(\xi) e^{-ix \xi} d\xi\]
for a smooth and rapidly decaying function $\phi$. We have
\[\chi\left( \frac{\log m_j}{\log R}\right)=\int_{-\infty}^{\infty}\phi(\xi_j) m_j^{-z_j} d\xi\;, \quad\mathrm{where} \;\;\;z_j = \frac{1+i\xi_j}{\log R}.\]
Setting up the corresponding notation involving $m_j',z_j'$ and $\xi_j'$, (\ref{eq:prev2}) becomes
\begin{align}\label{eq:prev3}
\left(\frac{\phi(W)}{W} \log R\right)^{|J|}\sum_{j\in J}& \sum_{m_j, m_j'\geq 1}  \int_{-\infty}^\infty \dots  \int_{-\infty}^\infty
\left( \prod_{j \in J}\mu(m_j)\mu(m_j') m_j^{-z_j}m_j'^{-z_j'} \phi(\xi_j)\phi(\xi_j')d\xi_j d\xi_j'\right)\\
&\prod_{p\leq R^{\log R}} c_p((WP_j+b_1)_{j \in J_1, r_j(p)=1},(WP_j+b_2)_{j \in J_2, r_j(p)=1})\nonumber,
\end{align}
which can be rewritten in the form
\[\left(\frac{\phi(W)}{W} \log R\right)^{|J|} \int_{-\infty}^\infty \dots  \int_{-\infty}^\infty \prod_{p\leq R^{\log R}}  E_p
\left( \prod_{i \in J}\phi(\xi_j)\phi(\xi_j')d\xi_j d\xi_j'\right),\]
where
\[E_p=\sum_{j\in J} \sum_{m_j, m_j' \in \{1,p\}} \left( \prod_{j \in J}\mu(m_j)\mu(m_j') m_j^{-z_j}m_j'^{-z_j'}\right)c_p((WP_j+b_1)_{j \in J_1, r_j(p)=1},(WP_j+b_2)_{j \in J_2, r_j(p)=1}).\]
One now approximates the Euler factor $E_p$ by 
\[E_p'= \prod_{j \in J} \frac{(1-p^{-(1+z_j)})(1-p^{-(1+z_j')})}{1-p^{-(1+z_j+z_j')}}\]
using a series of claims for different types of primes $p$, whose proofs we shall postpone until the end of the section.
%:claim:small
\begin{claim}[Small primes]\label{claim:small}
\[\prod_{p<w} \frac{E_p}{E_p'}=\left(\frac{W}{\phi(W)} \right)^{|J|}\left(1+o(1)\right)\]
\end{claim}
%:claim:bad
\begin{claim}[Bad but not terrible primes]\label{claim:bad}
\[\prod_{\stackrel{w<p\leq R^{\log R}}{p \;\mathrm{bad \;but \;not \;terrible}}} \frac{E_p}{E_p'}=1+O \left( \textup{Exp} \left( O \left( \sum_{p \in \mathcal{P}_b} p^{-1} \right) \right) \right)\]
\end{claim}
%:claim:good
\begin{claim}[Good primes]\label{claim:good}
\[\prod_{\stackrel{w<p\leq R^{\log R}}{p \;\mathrm{good}}} \frac{E_p}{E_p'}=1+o(1)\]
\end{claim}

Together Claims \ref{claim:small}, \ref{claim:bad} and \ref{claim:good} imply that if there are no terrible primes $>w$, then
\[ \prod_{w< p\leq R^{\log R}} \frac{E_p}{E_p'}=\left(\frac{W}{\phi(W)} \right)^{|J|}\left(1+o(1)+O \left( \textup{Exp} \left( O \left( \sum_{p \in \mathcal{P}_b} p^{-1} \right) \right) \right)\right).\]
The proof of Proposition \ref{prop:corr} is now completed, exactly as in \cite[Proposition 10.1]{tao-ziegler}, using some elementary theory of the Riemann $\zeta$ function, as well as the rapid decay of $\phi$.
\end{proofof}

Finally, for completeness, we give the proofs of  Claims \ref{claim:small}, \ref{claim:bad} and \ref{claim:good}. They rely on the rather elementary \cite[Lemma 9.5]{tao-ziegler}, which itself is proved with the help of a combinatorial Nullstellensatz \cite[Appendix D]{tao-ziegler}. 
%:proofof:claim:small
\begin{proofof}{of Claim \ref{claim:small}}
If $p<w$, then for all $j \in  J_i$ we have $WP_j+b_i \equiv b_i \not\equiv 0 (p)$, so for such $p$ the local factor $c_p((WP_j+b_1)_{j \in J_1, r_j(p)=1},(WP_j+b_2)_{j \in J_2, r_j(p)=1})$ is equal to 0 unless the set $\{j \in J: r_j(p)=1\}$ is the empty set, in which case it equals 1. The former case happens if and only if all $m_j$ are equal to 1, so that $E_p=1$.
A direct computation, using the fact that $w$ goes to infinity much more slowly than $R$, shows that $E_p'=(1-p^{-1})^{|J|}+o(1)$. The estimate for the product over all primes $p<w$ of $E_p/E_p'$ then follows from the fact that $\prod_{p<w}(1-p^{-1})^{-1}=W/\phi(W)$.
\end{proofof}
%:proofof:claim:bad
\begin{proofof}{of Claim \ref{claim:bad}}
If $p>w$ is bad but not terrible, then \cite[Lemma 9.5 (b)]{tao-ziegler} implies that the local factor $c_p$ of any non-trivial family is $O(p^{-1})$. Hence the sum defining $E_p$ has a contribution of 1 from $m_j=m_j'=1$ for all $j \in J$, and a contribution of $O(p^{-1})$ from all other terms, so that $E_p=1+O(p^{-1})$. Also, we find by Taylor expanding that
\begin{equation}\label{eq:epprime}
E_p'=1-\frac{1}{p} \sum_{j \in [J]} \left( \frac{1}{p^{z_j}}+\frac{1}{p^{z_j'}} \right)+\frac{1}{p} \sum_{j \in [J]} \frac{1}{p^{z_j+z_j'}} + O\left(\frac{1}{p^2}\right)=1+O\left(\frac{1}{p}\right),
\end{equation}
where the last inequality follows from the fact that $z_j, z_j'$ have real part $1/\log R >0$.
It follows that 
\[\prod_{\stackrel{w<p\leq R^{\log R}}{p \;\mathrm{bad \;but \;not \;terrible}}} \frac{E_p}{E_p'}= \prod_{\stackrel{w<p\leq R^{\log R}}{p \;\mathrm{bad \;but \;not \;terrible}}}\left(1+O\left(\frac{1}{p}\right)\right) = \prod_{\stackrel{w<p\leq R^{\log R}}{p \;\mathrm{bad \;but \;not \;terrible}}}\exp\left(O\left(\frac{1}{p}\right)\right),\]
which equals
\[\exp\left( O \left( \sum_{p \in \mathcal{P}_b} p^{-1} \right)   \right) = 1+O \left( \textup{Exp} \left( O \left( \sum_{p \in \mathcal{P}_b} p^{-1} \right) \right) \right).\]
\end{proofof}
%:proofof:claim:good
\begin{proofof}{of Claim \ref{claim:good}}
If $p>w$ is good, then by  \cite[Lemma 9.5 (c)]{tao-ziegler} the local factor $c_p$ of a non-trivial family consisting of precisely 1 polynomial is $p^{-1}+O(p^{-2})$. Such a non-trivial family arises in the case when $m_j=1, m_j'=p$ or  $m_j=p, m_j'=1$ or $m_j=p, m_j'=p$ for exactly one $j \in [J]$. In all remaining cases, $c_p=O(p^{-2})$ by \cite[Lemma 9.5 (d)]{tao-ziegler}. This implies that
\[E_p=1-\left(\frac{1}{p}+O\left(\frac{1}{p^2}\right)\right)  \sum_{j \in [J]} \left(\frac{1}{p^{z_j}}+\frac{1}{p^{z_j'}}\right) + \left(\frac{1}{p}+O\left(\frac{1}{p^2}\right)\right) \sum_{j \in [J]} \frac{1}{p^{z_j+z_j'}} + O\left(\frac{1}{p^2}\right), \]
so that by (\ref{eq:epprime}) for good primes $p$,
\[\frac{E_p}{E_p'}= 1+ O\left(\frac{1}{p^2}\right).\]
The statement in Claim \ref{claim:good} now follows from the fact that the product $\prod_p(1+O(p^{-2}))$ is convergent and $w$ tends to infinity.
\end{proofof}

\end{document}